\documentclass[10pt,reqno]{amsart}
\usepackage{amssymb,amsmath,amsthm,mathrsfs}
\usepackage{mathtools}
\usepackage{amscd}
\usepackage{subcaption}
\usepackage{graphicx}
\usepackage{color}
\usepackage[dvipsnames]{xcolor}
\numberwithin{equation}{section}
\usepackage{enumitem}
\usepackage[pagebackref=true, colorlinks=true, citecolor=blue]{hyperref}
\usepackage[misc]{ifsym}
\usepackage{epsfig} 
\usepackage{epstopdf} 
\usepackage[colorlinks=true]{hyperref}
\usepackage{mathabx}
\usepackage[normalem]{ulem}
\usepackage{mathtools}
\usepackage[tableposition=top]{caption}
\usepackage{booktabs,dcolumn}
\usepackage{stmaryrd}
\usepackage{cancel}
\usepackage{ulem}

\DeclareFontFamily{OT1}{rsfs}{}
\DeclareFontShape{OT1}{rsfs}{n}{it}{<-> rsfs10}{}
\DeclareMathAlphabet{\mathscr}{OT1}{rsfs}{n}{it}

\theoremstyle{plain}

\newtheorem{theorem}{Theorem}[section]
\newtheorem{proposition}[theorem]{Proposition}
\newtheorem{lemma}{Lemma}[section]

\theoremstyle{definition}

\newtheorem{definition}[theorem]{Definition}
\newtheorem{remark}{Remark}[section]

\def\be{\begin{equation}}
\def\ee{\end{equation}}

\newcommand{\bx}{\mathbf{x}}

\def\bes{\begin{equation*}}
\def\ees{\end{equation*}}
\def\bea{\begin{equation} \begin{aligned}}
\def\eea{\end{aligned} \end{equation}}
\def\beas{\begin{equation*} \begin{aligned}}
\def\eeas{\end{aligned} \end{equation*}}

\def\bi{\begin{itemize}}
\def\ei{\end{itemize}}

\newcommand{\ignore}[1]{}

\newcommand\eps{\varepsilon}

\newcommand{\pa}{\partial}
\newcommand{\na}{\nabla}

\newcommand\veps{\varepsilon}

\newcommand\mcL{\mathcal{L}}

\newcommand{\stkout}[1]{\ifmmode\text{\sout{\ensuremath{#1}}}\else\sout{#1}\fi}
\usepackage{algorithm, algorithmic}
 \let\oldalign\align
 \let\oldendalign\endalign
 \renewenvironment{align}
 {\linenomathNonumbers\oldalign}
 {\oldendalign\endlinenomath}

 \let\oldgather \gather
 \let\oldendgather\endgather
 \renewenvironment{gather}
 {\linenomathNonumbers\oldgather}
 {\oldendgather\endlinenomath}

\title[Sharp Material Interface Limit of the Darcy-Boussinesq System]{Sharp Material Interface Limit of the Darcy-Boussinesq System} 

\author{Hongjie Dong $^{{\href{mailto:hongjie_dong@brown.edu}{\textrm{\Letter}}}3}$}
\address{$^3$ Division of Applied Mathematics, Brown University, 182 George Street, Providence, RI 02912, USA}

\author{Xiaoming Wang $^{{\href{mailto:wxm.math@outlook.com}{\textrm{\Letter}}}1,2}$}
\address{$^1$ School of Mathematical Sciences, Eastern Institute of Technology, Ningbo, China}
\address{$^2$ Department of Mathematics and Statistics, Missouri University of Science and Technology, Rolla, MO, USA
}

\thanks{This work is supported by the NSF under agreement DMS2350129 (Dong) and NSFC 12271237 (Wang)}

\thanks{$^*$Corresponding author: Xiaoming Wang}

\begin{document}


 


\begin{abstract}
We investigate the sharp material interface limit of the Darcy-Boussinesq model for convection in layered porous media with diffused material interfaces, which allow a gradual transition of material parameters between different layers.
We demonstrate that as the thickness of these transition layers approaches zero, the conventional sharp interface model with interfacial boundary conditions, commonly adopted by the fluids community, is recovered under the assumption of constant porosity.
Our results validate the widely used sharp interface model by bridging it with the more physically realistic case of diffused material interfaces.
This limiting process is singular and involves a boundary layer in the velocity field.
Our analysis requires delicate estimates for elliptic and parabolic equations with discontinuous coefficients,  and the subtle validation of the boundary layer.
\end{abstract}
 \maketitle

\section{Introduction}


Convection in porous media, which is highly relevant for geophysical applications and many engineering problems, has been the focus of numerous researchers; see, for instance, Nield and Bejan \cite{nield2017convection}. Several mathematical works also address this area \cite{doering1998jfm, fabrie1986aam, fabrie1996aa, ly1999jns, oliver2000gevrey, otero2004jfm, CNW2024}. However, all these works—except for a few sections in the book by Nield and Bejan and \cite{CNW2024}—deal with homogeneous porous media.

On the other hand, many natural and engineered porous media are non-homogeneous. Such inhomogeneity may arise in layered structures, where the domain consists of multiple layers with permeability and other physical parameters that vary from layer to layer. These layered structures can be formed naturally due to processes such as sedimentary deposition, volcanic and tectonic activity, and groundwater flow and dissolution, among others.

Layered porous media are also closely related to the underground carbon dioxide ($CO_2$) sequestration technology, which has attracted significant attention from the fluid dynamics and environmental science communities \cite{bickle2007modelling, huppert2014ar, hewitt2014jfm, hewitt2020jfm, hewitt2022jfm, mckibbin1980jfm, mckibbin1981heat, mckibbin1983thermal, sahu2017tansp, salibindla2018jfm, wooding1997convection}. The review paper by Huppert and Neufeld \cite{huppert2014ar} provides an excellent survey of this topic. See \cite{BS2024} and the references therein for recent work on the utility of layered porous media in enhancing heat transport.

A simple treatment of a layered domain is to assume that the material parameters remain constant within each layer. These parameters undergo a sudden jump at the interfaces between adjacent layers, and the models in neighboring layers are coupled through certain ad-hoc interfacial boundary conditions; see \cite{hewitt2014jfm, hewitt2020jfm, hewitt2022jfm, mckibbin1980jfm, mckibbin1981heat, mckibbin1983thermal, CNW2024} and \eqref{IC}.

Such a “sharp” material interface approach is an idealization, considering the ever-present microscopic diffusion processes. Additionally, although the homogeneous model has been derived rigorously via homogenization method, the mathematical derivation of the layered model is not known. Moreover, some of the interfacial boundary conditions associated with the sharp material interface model appear somewhat mysterious. For instance, only the continuity of the velocity normal to the interface is enforced, whereas the tangential velocity is allowed to be discontinuous—an observation that is counterintuitive for a slow-evolving fluid system.

On the other hand, a physically more realistic model would account for diffusion and allow for a narrow transition region where the material parameters gradually change from a constant value in one layer to another constant value in an adjacent layer.
Following \cite{diffuse98, SMCTSH17}, we refer to this model as the diffuse material interface model. A fundamental physical question, then, is whether the sharp interface model, together with the interfacial boundary conditions, can be recovered as the limiting case of the diffuse material interface model when the thickness of the transition region approaches zero.

The main question we address in this manuscript is whether the solutions of the diffused interface model \eqref{2.1}-\eqref{2.3} converge to the solutions of the sharp interface model \eqref{1.1}-\eqref{1.3} as the width of the diffuse interface approaches zero, i.e., $\veps\rightarrow 0$, together with the appropriate interfacial boundary conditions.
Our main accomplishment is providing an affirmative answer to this question, thereby establishing rigorous mathematical and physical foundations for the sharp interface model along with its interfacial boundary conditions.

We also discover that the sharp material interface limit is singular, in the sense that the velocity field develops a boundary layer-like structure around each interface when the transition layer width is small. The singular nature of the problem makes the proof non-trivial.
We leverage the anisotropic regularity of the solutions—particularly their regularity in the direction(s) parallel to the interface—for both the sharp material interface model, which has discontinuous material parameters, and the diffuse interface model, whose material properties exhibit large gradients in the direction perpendicular to the interfaces.
In particular, the validity of the boundary layer expansion is established using a novel embedding result for function spaces with anisotropic fractional regularity.

The rest of the paper is organized as follows. In Section 2, we introduce both the sharp and diffuse interface models. Section 3 presents some preliminary estimates for both models. Section 4 is devoted to the convergence analysis of the concentration and pressure fields. The boundary layer structure of the velocity field is discussed in Section 5. Concluding remarks are offered in Section 6. The novel embedding result for function spaces with anisotropic fractional regularity is included in the Appendix.

\section{The diffuse and sharp material interface models}

We present the mathematical models together with the associated boundary and interfacial boundary conditions in this section. We also define the concept of weak solution to these models.
A novel natural space associated with the solutions to the sharp material interface model is introduced, and its basic properties such as piecewise smoothness and Sobolev-type inequalities are also presented. 
For simplicity, we consider the idealized 2-d domain $\Omega=(0,L) \times(-H, 0)$
\subsection{The sharp interface model}
We first introduce the idealized `sharp' interface model with $\ell-1$ interfaces located at $z=z_{j}\in (-H, 0), j=1, \cdots,\ell-1$. We also denote $z_0=0, z_\ell=-H$. These interfaces divide the domain into $\ell$ `layers' or `strips, ' denoted $\Omega_j$:
\begin{equation}
\Omega_j=\left\{\textbf{x}=(x,z) \in \Omega \mid z_{j}<z<z_{j-1}\right\}, \quad 1 \leq j \leq   \ell .\label{domain}
\end{equation}
For convection in this layered domain, the well-accepted ad-hoc governing equations in $\Omega$ are the following Darcy-Boussinesq system (with the usual Boussinesq approximation) \cite{nield2017convection}:
\begin{eqnarray}
\operatorname{div}(\mathbf{u})=0, \quad \mathbf{u}=\left(u_1, u_2\right), \label{1.1}\\
\mathbf{u}=-\frac{K}{\mu}\left(\nabla P+\rho_0(1+\alpha \phi) g \mathbf{e}_z \right),\label{1.2} \\
b \frac{\partial \phi}{\partial t}+\mathbf{u} \cdot \nabla \phi-\operatorname{div}(b D \nabla \phi)=0.\label{1.3}
\end{eqnarray}
Here  $\mathbf{u}, \phi$ and $P$ are the unknown fluid velocity, concentration, and pressure, respectively; $\rho_0, \alpha,  \mu, g$ are the constant reference fluid density,   constant expansion coefficient,  constant dynamic viscosity, and the gravity acceleration constant, respectively; and $\mathbf{e}_z $ stands for the unit vector in the $z$ direction. In addition,
$K, b, D$ represent the permeability, porosity, and diffusivity coefficients respectively, which are assumed to be constant within each strip/layer  $\Omega_j$, i.e.,
\begin{equation}
K=K(\bx)=K_j,\quad b=b(\bx)=b_j,\quad D=D(\bx)=D_j, \quad \bx \in \Omega_j, \quad 1 \leq j\leq \ell,   \label{coefficients}
\end{equation}
for a set of positive constants $\left\{K_j, b_j, D_j \right\}_{j=1}^{\ell}$. Here $K$ and $D$ could be positive symmetric tensors in the more general case.

System \eqref{1.1}-\eqref{1.3} models convection in a layered porous media where each layer is of different permeability/porosity/diffusivity. On the interfaces $z= z_j$,  we assume
\begin{eqnarray}\label{IBC_up}
    \mathbf{u}\cdot \mathbf{e}_z ,\,\, P, \,\, \phi  \text{ are continuous at } z= z_j, \, 1 \leq j \leq  \ell-1 .
  \label{IBC}
\end{eqnarray}
These conditions are needed in order to ensure conservation of fluid concentration mass, and the balance of normal force. Notice that we do not impose the continuity of the tangential velocity.
Interfacial boundary condition \eqref{IBC} together with the concentration transport and diffusion equation \eqref{1.3} imply that the solution $\phi$ cannot be smooth over the whole domain in general unless $bD$ is a constant \cite{mckibbin1980jfm, mckibbin1981heat}.
{ Moreover, \eqref{1.3} also implies, for solution that is piecewise smooth,}
\be
  bD\frac{\phi}{\partial z} \ \mbox{continuous\ at}\  a=z_j, j=1, \cdots,\ell-1.
  \label{IBC-2}
\ee
{ This interfacial boundary condition will be implicit in the weak formulation that we introduce later.}

The system \eqref{1.1}-\eqref{1.3} can be reformulated as a coupled system of subsystems on each layer/strip, with \eqref{IBC} and \eqref{IBC-2} serving as the coupling conditions. This is reminiscent of classical transmission problems, but with the added challenge of nonlinearity as well as discontinuous thermal diffusivity. See, for instance, \cite{Lions56, St56, LU68, LN03}.

System \eqref{1.1}-\eqref{1.3} is supplemented with the initial condition
\be\label{IC}
\phi( x,z , 0)=\phi_0(x,z),
\ee
and the boundary conditions (prescribed concentration)
\begin{equation*}
\begin{split}
\mathbf{u}\cdot e_z \left(x, 0 ; t\right)=\mathbf{u}\cdot e_z \left(x,-H ; t\right)=0, \\
\phi\mid_{z=0}=C_0, \quad \phi\mid_{z=-H}=C_1
\end{split}
\end{equation*}
together with periodicity in the horizontal direction(s) $ x$ ($x=(x_1,x_2)$ in the three-dimensional case).

Let $\phi_b(z)$ be a smooth function satisfying  $\phi_b(0)=C_0, \phi_b(-H)=C_1$, and $D\frac{d\phi_b}{dz}=0, z\in (\frac{-H+ z_{\ell-1}}{2}, \frac{z_1}{2})$.
    Introducing $\widetilde{\phi}(x,z, t)=\phi(x,z, t)+\phi_b(z)$, the boundary condition for $\phi$ can be homogenized. The change of variable also introduces extra lower order terms involving $\phi_b$  in \eqref{1.3}. These lower order  terms do not impose any genuine difficulty in our analysis. Hereafter for simplicity,   and without loss of generality, we assume that $C_0=C_1=0$ . Thus,
the boundary condition becomes
    \begin{equation}\label{hBC}
       \phi\mid_{z=0}=0,\,\,   \phi\mid_{z=-H}=0; \quad \mathbf{u}\cdot e_z \left(x, 0 ; t\right)=\mathbf{u}\cdot e_z \left(x,-H ; t\right)=0.
    \end{equation}

By setting $\widetilde{P}=P-\rho_0  g z$, and omitting $\sim$ for simplicity, we
may rewrite \eqref{1.2} as
\begin{equation}
\mathbf{u}=-\frac{K}{\mu}\left(\nabla P+ \alpha \rho_0 g \phi \mathbf{e}_z \right). 
\end{equation}
Notice that the interfacial conditions \eqref{IBC_up} remain unchanged under the change of variable.
We will adopt the new formulation  hereafter. Moreover, we will assume $\mu=1, \alpha\rho_0g=1, b=b^\veps=1$ for convenience. 
Hence, the simplified transformed {\bf sharp interface model} takes the form
\begin{eqnarray}
\operatorname{div}(\mathbf{u})=0, \quad \mathbf{u}=\left(u_1, u_2\right), \label{1.1'}\\
\mathbf{u}=-K\left(\nabla P+ \phi \mathbf{e}_z \right),\label{1.2'} \\
 \frac{\partial \phi}{\partial t}+\mathbf{u} \cdot \nabla \phi-\operatorname{div}( D \nabla \phi)=0,\label{1.3'}
\end{eqnarray}
together with the interfacial boundary conditions \eqref{IBC}, \eqref{IBC-2}, boundary condition \eqref{hBC}, and initial condition  \eqref{IC}.

\subsection{The diffuse interface model}
While the idealized model allows the physical parameters make sudden changes at $z=z_j, 1\le j\le \ell-1$, a physically more realistic setting would allow the material parameters make gradual changes, albeit over a thin region, with thickness $\veps\ll 1$, due to local diffusion. We assume that
$$\veps \ll z_j-z_{j+1}, \quad \forall j=0, \cdots,\ell-1.$$
This implies that the gradual transition layers are much thinner than the macroscopic layers. 

The governing equations then take the form 
\begin{eqnarray}
\operatorname{div}(\mathbf{u}^\veps)=0, \quad \mathbf{u}^\veps=\left(u^\veps_1, u^\veps_2\right), \label{2.1}\\
\mathbf{u}^\veps=-\frac{K^\veps}{\mu}\left(\nabla P^\veps+\rho_0(1+\alpha \phi^\veps) g \mathbf{e}_z \right),\label{2.2} \\
b^\veps \frac{\partial \phi^\veps}{\partial t}+\mathbf{u}^\veps \cdot \nabla \phi^\veps-\operatorname{div}(b^\veps D^\veps \nabla \phi^\veps)=0.\label{2.3}
\end{eqnarray}
Here  $\mathbf{u}^\eps, \phi^\veps$ and $P^\veps$ are the unknown fluid velocity, concentration, and pressure for the diffuse interface model respectively.
$K^\veps, b^\veps, D^\veps$ represent the permeability, porosity, and diffusivity coefficients that are essentially constant in each of the $\Omega_j$, but with a narrow thin layer of width $2\veps$ that allows a smooth transition between the neighboring strips.
Namely,
\begin{equation}
        \label{diffuse_coefficients} 
    K^\veps=K_j,\quad b^\veps=b_j,\quad D^\veps=D_j, \quad z_j+\veps \le z \le z_{j-1}-\veps, \quad 1 \leq j\leq \ell;
\end{equation}
with $ K^\veps, b^\veps, D^\veps$ smooth over $\Omega.$    

For simplicity, we will assume $K^\veps, D^\veps$ be continuous piecewise linear functions in the $z$ variable, i.e.,
$$K^\veps(z)=K_{j-1}+(z-z_j+\veps)\frac{K_j-K_{j-1}}{e\veps}, z\in (z_j-\veps, z_j+\veps). $$
In addition, we assume  $b^\veps\equiv const$.
No interfacial conditions are needed since all variables are continuous in this diffusive model.

Utilizing the same homogenization of the boundary condition method as well as the change of variable procedure as for the sharp interface model, we deduce a simplified {\bf diffusive interface model} after ignoring the same kind of lower order terms:
\begin{eqnarray}
\operatorname{div}(\mathbf{u}^\veps)=0, \quad \mathbf{u}^\veps=\left(u^\veps_1,\cdots, u^\veps_d\right), \label{2.1'}\\
\mathbf{u}^\veps=-K^\veps\left(\nabla P^\veps+ \phi^\veps \mathbf{e}_z \right),\label{2.2'} \\
 \frac{\partial \phi^\veps}{\partial t}+\mathbf{u}^\veps \cdot \nabla \phi^\veps-\operatorname{div}( D^\veps \nabla \phi^\veps)=0,\label{2.3'}
\end{eqnarray}
together with the boundary condition \eqref{hBC} (with the variable replaced by the variable with superscript $\veps$) and initial condition  \eqref{IC} (with $\phi$ replaced by $\phi^\veps$).

\noindent{\bf Remark:} Note that the diffuse interface model terminology that we have adopted here is different from the conventional diffusive interface model where the so-called diffusive interface is the diffused approximation of an otherwise moving sharp interface/front associated with the system, such as the interface between two immiscible fluids 
\cite{diffuse98}.
The diffusive interface in this paper is the diffused version of an otherwise fixed sharp interface that separates two regions of the domain with different material parameters following the pioneering work of Hughes and collaborators \cite{SMCTSH17}.


\section{Preliminary Estimates}

\subsection{Weak formulation}
While the solution to the diffuse interface model  \eqref{2.1}-\eqref{2.3} can be classical, the solution to the sharp interface model \eqref{1.1}-\eqref{1.3} must be understood in the weak sense due to the discontinuity in the material parameters. This is similar to the study of transmission problems if formulated as global problems. See, for instance \cite{Lions56, St56, LU68, LN03}.

Let $L^p(\Omega)$ and $H^k(\Omega)$ denote the usual $L^p$-Lebesgue space of integrable functions and $H^k$ Sobolev spaces that are periodic in the horizontal direction(s), respectively, for $1 \leq p \leq \infty$ and $k \in \mathbb{R}$. The inner product in $L^2(\Omega)$ will be denoted by $(\cdot , \cdot)$. Let
 \begin{eqnarray*}
\mathcal{V}:=\left\{\phi \in C(\bar{\Omega}) \text{ and }\phi\mid_{\Omega_j}\in C^{\infty}(  \bar{\Omega}_j  ),\, \, 1\leq  j \leq \ell  \mid \phi \text{ is periodic in $x$ and satisfies }\eqref{IBC},   \eqref{hBC}\right\}, \\
V:=\text { Closure of } \mathcal{V} \text { in the } H^1-\text {norm},\\
H:=\text { Closure of } \mathcal{V} \text { in the } L^2-\text {norm},
 \end{eqnarray*}
and let us denote the $L^2$-norm of $H$ by $\|\cdot\|_H$, and the norm of $V$ by $\|\cdot\|_V$. The inner product of $H$ is exactly the inner product of $L^2(\Omega)$. Notice that due to the boundary conditions \eqref{hBC}, the Poincar\'{e} inequality implies that the $V$-norm and the $H^1$-Sobolev norm are equivalent and thus,  when combined with the lower and upper bounds on $b$ and $D$, we can define the $V$ norm in the following manner
 \begin{equation}\label{V-norm}
 \|\phi\|_V=\|\sqrt{bD}\nabla \phi\|_{L^2(\Omega)}.
 \end{equation}
We also recognize that $V=H^1_{0,per}(\Omega)$, the subspace of $H^1(\Omega)$ that vanishes at $z=0, -H$ and periodic in the horizontal direction(s).
We denote the dual space of $V$ by $V^* $ with norm $\|\cdot\|_{V^{*}}$. The symbol $\langle\cdot, \cdot\rangle$ will stand for the duality product between $V$ and $V^* $.

Let us also define
\begin{eqnarray*}
\tilde{\mathcal{V}}:=\left\{\mathbf{u} \in C(\bar{\Omega})^d \text{ and }\mathbf{u}\mid_{{\Omega}_j}\in C^{\infty}(\bar{\Omega}_j)^d,\, \, 1 \leq j \leq  \ell  \mid \mathbf{u}\text { satisfies }\eqref{1.1},\eqref{IBC_up}, \eqref{hBC}\right\} \\
\mathbf{H}:=\text { Closure of } \tilde{\mathcal{V}} \text { in the } L^2-\text {norm. }
\end{eqnarray*}
Utilizing the divergence free property of the velocity $\mathbf{u}$, the continuity of the pressure $P$, the pressure can be formulated as a function of the concentration $\phi$ through the following second order elliptic equation with discontinuous coefficient:
\begin{equation}\label{eqP}
\begin{cases}
-\text{div} \Big( {K}\nabla P\Big) =     \text{div} \Big(   K \phi \mathbf{e}_z  \Big) \quad \text{ in } \Omega,   \\
 \frac{\partial P}{\partial z }(x,0) =-\phi(x,0),\quad   \frac{\partial P}{\partial z }(x,-H)=-\phi(x,-H) .
\end{cases}
\end{equation}
Here the Neumann boundary condition is understood in the weak sense. The solvability of the equation in $H^1(\Omega)$ 
for any $\phi\in H$ is guaranteed by the Lax-Milgram theorem using the following weak formulation:
\begin{equation}
                        \label{eq1.02}
\int_\Omega{K} \nabla  P\cdot \nabla q d\bx =   \int_\Omega     K \phi \mathbf{e}_z \cdot \nabla q  d\bx .
\end{equation}
Uniqueness can be obtained if we restrict to the mean zero subspace of $H^1$.
For piecewise smooth $\phi$, i.e., smooth in each layer $\Omega_j$, one can show that $P$ is also piecewise smooth, and the homogeneous Neumann boundary can be verified provided that the concentration $\phi$ vanishes at the top and bottom of the domain in the trace sense. More specifically, if $\phi\in V$, then $P$ is piecewise $H^2$ and therefore, the boundary conditions in \eqref{eqP} become the homogeneous boundary conditions in the trace sense. 

\begin{lemma}
    \label{lem3.1}
(1) When $\phi\in L^p(\Omega)$ for $p\in (1,\infty)$, we have $\nabla P\in L^p(\Omega)$ and 
\begin{equation}
    \label{eq9.49}
\|\nabla P\|_{L^p(\Omega)}\le C(d,p,K_j)\|\phi\|_{L^p(\Omega)},
\end{equation}
where $C(d,p,K_j)\le C(d,K_j)p^\beta$ for $p\in (2,\infty)$ with some $\beta=\beta(d)>0$.

(2) When $\phi$ is piecewise $C^\alpha$ H\"older continuous for some $\alpha\in (0,1)$, $P$ is globally Lipschitz in $\Omega$ and piecewise $C^{1,\alpha}$. Furthermore, we have
\begin{equation*}
    \|\nabla P\|_{L^\infty(\Omega)}+[\nabla P]_{C^\alpha(\Omega_j)}\le C(d,\alpha,K_j)\sum_j\|\phi\|_{C^\alpha(\Omega_j)}
\end{equation*}
\end{lemma}
\begin{proof}
See \cite{D10}. We note that in view of the proof in \cite{D10} and the linear dependence on $p$ of the implicit constant in the Fefferman-Stein theorem (see \cite[$\S$3.2]{MR2435520} and \cite[$\S$IV.2]{MR1232192}, $C(d,p,K_j)$ in \eqref{eq9.49} has a polynomial growth with respect to $p\ge 2$.
For (2), see \cite{D12}.
\end{proof}

Define the bilinear forms $A(\phi,\psi )$  and $ B(\mathbf{u}, \phi,\psi )$  for some $\mathbf{u}\in \mathbf{H}\cap L^p(\Omega)^d$, with a $p>2$:
\begin{eqnarray*}
A(\phi,\psi )&=  \left( bD\nabla \phi, \nabla \psi\right), \quad \forall \phi, \psi \in V , \\
B(\mathbf{u}, \phi,\psi )&= \int_{\Omega}\left(\mathbf{u}\cdot \nabla \phi \right)\psi d\bx, \quad \forall \phi, \psi \in V. \notag
\end{eqnarray*}
The weak solution to the system \eqref{1.1'}-\eqref{1.3'} is defined as follows.
\begin{definition}\label{def1}
Let $\phi_0 \in L^2(\Omega)$  be given, and let $T>0$. A weak solution of \eqref{1.1'}-\eqref{1.3'}, subject to the boundary conditions  \eqref{IBC}, \eqref{hBC} together with the periodic conditions in the horizontal direction, and the initial condition \eqref{IC} on the interval $[0,T]$ is a triple $(\mathbf{u},\phi, P)$, satisfying
$$
\phi \in L^2(0, T ; V) \cap  L^{\infty}\left(0, T ; L^2(\Omega)\right) \text{ and }{\partial \phi}/{\partial t} \in L^2\left(0, T ; {V^* }\right),   \notag
$$
and
\begin{equation}\label{def-eq}
\left(\phi(t_2), \psi\right)-\left(\phi(t_1), \psi\right)+\int_{t_1}^{t_2}A( \phi, \psi) dt+\int_{t_1}^{t_2}B(\mathbf{u},  \phi,\psi) d t=0,\end{equation}
$\forall  \psi \in V, t_1, t_2 \in[0,T]$,  and  $\phi(0)=\phi_0 $ in $L^2(\Omega)$,
where $\mathbf{u} \in L^2\left(0, T ; L^p(\Omega)\right)$ and $P \in  L^2\left(0, T ; W^{1,p}(\Omega)\right) $, with $2 \le p <\infty$ if $d=2$, and $p=6$ if $d=3$,  are given by Darcy's law, i.e., \eqref{1.2} and \eqref{eqP}, respectively.

\end{definition}

\subsection{Estimates on  the sharp interface model}
While the weak solution to the sharp interface model can be derived via energy methods and the elliptic regularity for systems with partially small BMO coefficients \cite{dongli21jfa}, higher order estimates are needed for the vanishing material interface limit analysis.

For higher regularity of the solutions to \eqref{1.1'}-\eqref{1.3'}, we recall the following space associated with the 
 principal differential operator of the system \cite{CNW2024}, i.e.,
  $\mcL= -\text{div}(bD \na)$, subject to the interfacial boundary conditions \eqref{IBC},  and \eqref{hBC}  together with the periodic conditions in the horizontal direction.
 \begin{definition}\label{natural-space}
 Define
  \begin{equation}\label{W}
  W= \{\varphi\in V:  \pa_{x} \varphi \in H^1(\Omega), ~ bD \pa_{z} \varphi  \in H^1(\Omega) \},
 \end{equation}
endowed  with norm
\begin{equation}\label{W-norm}
 \|\varphi\|^2_{W}= \|\varphi \|^2_V+ \|\pa_{x}\varphi \|^2_{H^1(\Omega)} + \|bD\pa_{z}\varphi \|^2_{H^1(\Omega)}.
 \end{equation}
 \end{definition}

 \begin{remark}
                \label{rem3.1}
 $W$ is different from the classical $H^2$ space in general. Indeed, it is easy to verify the following equality.
 $$W\bigcap H^2=\{\phi\big| \phi\in H^2(\Omega),  \partial_z\phi|_{z=z_j}=0\, \, j=1, \cdots, l-1\}$$
 if $bD$ contains a jump discontinuity at $z=z_j, j=1, \cdots, l-1$.
 Moreover, if we make a change of variables in $z$ so that $d z/d \tilde z=bD$, then $\varphi\in W$ if and only if it is in $H^2$ with respect to the new variables $(x,\tilde z)$. Notice that this change of variable is Lipschitz but not smooth, unless $bD$ is smooth.
 \end{remark}
 
We recall the following lemma about an equivalent norm on $W$ associated with the operator $\mcL$ \cite{CNW2024} 
\begin{lemma}\label{remark-norm-equi}
There exist $C_l, C_u>0$ such that
\begin{equation}\label{remark-norm-equi-00}
C_l \|\mathcal{L} \phi\|_{L^2(\Omega)}\le \|\phi\|_{W} \leq C_u \|\mathcal{L} \phi\|_{L^2(\Omega)}, \quad \forall\phi\in W.
\end{equation}
Hence, the norm $\|\phi\|_{W} $ is equivalent to  $\|\mathcal{L} \phi\|_{L^2(\Omega)}.$

\end{lemma}

We also recall that the space $W$ enjoys some Sobolev type inequalities similar to the $H^2$ space \cite{CNW2024}. 
\begin{lemma} \label{remark-imbedding}
 $W$ is similar to $H^2$ in the sense that the following inequalities hold
 \begin{equation} \label{remark-imdedding-1}
   \| \na \varphi \|_{L^p(\Omega)}   \leq C  \|\varphi\|_{W},\, 2\leq  p<\infty,  \text{ for }  d=2, \,\,\,\text{ and } \,\,\,
    \| \na \varphi \|_{L^6(\Omega)}   \leq C  \|\varphi\|_{W}  \text{ for }  d=3.
    \end{equation}
and
 \begin{equation}\label{phi-infinity}
  \|   \varphi \|_{L^\infty(\Omega)}   \leq C  \|\varphi\|^\frac12\|\varphi\|_{W}^\frac12 \quad \text{ for }  d=2,
  \quad\text{ and}\quad
  \|   \varphi \|_{L^\infty(\Omega)}   \leq C  \|\varphi\|^\frac14\|\varphi\|_{W}^\frac34 \quad \text{ for }  d=3.
  \end{equation}
  \end{lemma}

\begin{proposition}
                \label{prop3.3}
For each $\phi_0\in H$, there exists a unique weak solution $(\textbf{u}, \phi, P)$  to problem \eqref{1.1'}-\eqref{1.3'} subject to the boundary conditions \eqref{hBC}, interfacial boundary condition \eqref{IBC}, together with the periodic conditions in the horizontal direction, and the initial condition. In addition,  we have the following regularity result.
\begin{enumerate}
    \item 
   If $\phi_0 \in V$,   we have, for any $T>0$,
   \begin{eqnarray}           \label{eq12.43}
\phi \in L^\infty(0,T; V ) \cap L^2(0,T; W),\quad \partial_t\phi\in L^2(0,T; L^2(\Omega)),
\\
                \label{eq10.12}
P\in L^\infty_t(W^{1,6}(\Omega)),\quad \frac{\partial P}{\partial t}\in L^2_t(H^1),\quad P\big|_{\Omega_j}\in L^\infty_t(H^2(\Omega_j))\cap L^2_t(H^3(\Omega_j))\, \forall j, 
\\
            \label{eq10.13}
\frac{\partial\mathbf{u}}{\partial t}\in L^2_t(L^2),\quad \mathbf{u}\big|_{\Omega_j}\in L^\infty_t(H^1(\Omega_j))\cap L^2_t(H^2(\Omega_j))\,\, \forall j.
    \end{eqnarray}
\item   If $\phi_0 \in W$, we have, for any $T>0$, 
    \begin{eqnarray}
        \label{eq10.30}
    \phi \in L^\infty_t(W),\quad  \partial_t \phi\in L^\infty_t(H) \cap L^2_t(V),
    \\
        \label{eq10.31}
  \nabla P\big|_{\Omega_j}\in L^\infty_t(H^2(\Omega_j)),
  \quad   \nabla \partial_t P\big|_{\Omega_j}\in L^\infty_t(L^2(\Omega))\cap  L^2_t(H^1(\Omega_j)),\\
        \label{eq10.32}
 \frac{\partial\mathbf{u}}{\partial t}\in L^\infty_t(L^2(\Omega))\cap  L^2_t(H^1(\Omega_j)), \quad \mathbf{u}\big|_{\Omega_j}\in L^\infty_t( H^2(\Omega_j)),\\
        \label{eq10.30b}
        \nabla \phi\in L^2_t(L^\infty(\Omega)).
  \end{eqnarray}
 \item If $\phi_0, \frac{\partial \phi_0}{\partial x}\in W$, we have, for any $T>0$,
\begin{eqnarray}
                    \label{eq10.33}
  \nabla \phi \in L^\infty_t(L^\infty(\Omega)), \quad \frac{\partial \phi}{\partial x} \in L^\infty_t(W), \quad \frac{\partial^2 \phi}{\partial t\partial x} \in L^\infty_t(L^2(\Omega))\cap L^2_t(V),
  \\
                \label{eq10.34}
  \nabla\frac{\partial P}{\partial x}\big|_{\Omega_j} \in L^\infty_t(H^2(\Omega_j)),\quad \nabla \partial_x\partial_t P\big|_{\Omega_j}\in L^\infty_t(L^2(\Omega))\cap  L^2_t(H^1(\Omega_j)),
  \\
                    \label{eq10.35}
  \frac{\partial^2\mathbf{u}}{\partial x\partial t}\in L^\infty_t(L^2(\Omega)), \quad \frac{\partial^2\mathbf{u}}{\partial x\partial t}\big|_{\Omega_j}\in L^2_t(H^1(\Omega_j)),\quad \frac{\partial\mathbf{u}}{\partial x}\big|_{\Omega_j}\in L^\infty_t( H^2(\Omega_j)),\\
            \label{eq10.33b}
\partial_x\nabla \phi\in L^2_t(L^\infty(\Omega)).
\end{eqnarray}
  \item Let $k\ge 2$. If $\phi_0, \frac{\partial \phi_0}{\partial x},  \frac{\partial^2 \phi_0}{\partial x^2},\ldots, \frac{\partial^k \phi_0}{\partial x^k}\in W$, we have, for any $T>0$,
\begin{eqnarray}
           \label{eq11.45}
\nabla\partial^{k-1}_x \phi \in L^\infty_t(L^\infty(\Omega)), \quad   \partial^k_x \phi \in L^\infty_t(W), \quad \partial_t\partial_x^k \phi \in L^\infty_t(L^2(\Omega))\cap L^2_t(V),
  \\
            \label{eq11.46}
  \nabla \partial_x^k P\big|_{\Omega_j} \in L^\infty_t(H^2(\Omega_j)),\quad \nabla \partial^k_x\partial_t P\big|_{\Omega_j}\in L^\infty_t(L^2(\Omega))\cap  L^2_t(H^1(\Omega_j)),
  \\
             \label{eq11.47}
  \partial_t\partial^k_x \mathbf{u}\in L^\infty_t(L^2(\Omega)),
  \quad \partial_t\partial^k_x \mathbf{u}\big|_{\Omega_j}\in L^2_t(H^1(\Omega_j)),
  \quad \partial^k_x\mathbf{u}\big|_{\Omega_j}\in L^\infty_t( H^2(\Omega_j)),\\
        \label{eq11.48}
       \partial_x^k\nabla \phi\in L_t^2(L^\infty(\Omega)).
\end{eqnarray}
\end{enumerate}
\end{proposition}
\begin{proof}
(1) The proof of \eqref{eq12.43} can be found in \cite{CNW2024}.
Now it follows from \eqref{eq12.43}, the Sobolev embedding $H^1\subset L^6$, and \eqref{eq9.49}, we have $P\in L^\infty_t(W^{1,6}(\Omega))$ and $\partial_t P\in L_t^2(H^1)$. By differentiating \eqref{eqP} in $x$ and using \eqref{eq9.49}, we get
$\partial_x P\in L^\infty_t(H^1(\Omega))$ as well as $P\in L^\infty_t(H^2(\Omega_j))$ for any $j$. Differentiating \eqref{eqP} twice in $x$ and using $\partial_x^2 \phi \in L^2_t(L^2(\Omega))$, we then obtain 
\begin{equation}
    \label{eq10.10}
\partial^2_x P\in L^2_t(H^1(\Omega)).
\end{equation}
Since for almost every $t$, $P$ satisfies $-\Delta P=\partial_z \phi$ in each subdomain $\Omega_j$, so $\partial_z^2 P=-\Delta_x P-\partial_z \Phi$ in $\Omega_j$. This and \eqref{eq10.10} imply $P\in L_t^2(H^3(\Omega_j))$ for each $j$. This completes the proof of \eqref{eq10.12}. Finally, \eqref{eq10.13} follows from \eqref{eq12.43}, \eqref{eq10.12}, and \eqref{1.2'}.


(2) Now we proceed with the proof of (2) on $\phi$.
We differentiate \eqref{1.1'}-\eqref{1.3'} in time, and we deduce
$$ \frac{\partial^2 \phi}{\partial t^2}+\mathbf{u} \cdot \nabla \frac{\partial\phi}{\partial t}
+\frac{\partial\mathbf{u}}{\partial t} \cdot \nabla \phi
-\operatorname{div}( D \nabla \frac{\partial\phi}{\partial t})=0, \quad\frac{\partial \phi}{\partial t}\big|_{z=0, -H}=0, \quad \frac{\partial\phi}{\partial t}\big|_{t=0}=-\mcL(\phi_0)-\mathbf{u}_0\cdot\nabla\phi_0.$$
By Lemma \ref{remark-norm-equi}, we have $\mcL(\phi_0)\in H$. By Lemmas \ref{remark-imbedding} and \ref{lem3.1}, we also have $\mathbf{u}_0\cdot\nabla\phi_0\in H$. 
Since $\phi_0\in W\subset L^\infty(\Omega)$, by the maximum principle, $\phi\in L^\infty(0,T;L^\infty)$.
Multiplying this equation by $\partial_t\phi$, integrating over the domain, and integrating by parts, 
utilizing the estimate 
$|\int_\Omega\partial_t\mathbf{u}\cdot\nabla\phi\partial_t\phi|\le \|\partial_t\mathbf{u}\|_{L^2}\|\nabla\partial_t\phi\|_{L^2}\|\phi\|_{L^\infty}$ and the maximum principle on $\phi$,
we deduce that $\partial_t\phi\in L^\infty(0,T; H)\cap L^2(0,T;V)$. 

The first point 
also implies $\mathbf{u}\cdot\nabla\phi\in L^\infty(0,T; L^{3/2}).$
Equation \eqref{1.3'} together with the estimate $\frac{\partial \phi}{\partial t} \in L^\infty(0,T; H)$ implies $\phi,\partial_x \phi, D\partial_z\phi\in L^\infty(0,T; W^{1,3/2}(\Omega))$ via the elliptic regularity theory. This further implies that $\mathbf{u}\cdot\nabla\phi\in L^\infty(0,T; L^2)$, which leads to the desired $L^\infty_t(W)$ estimate on $\phi$. 
This formal analysis can be rigorously justified through a Galerkin approximation using eigenfunctions of $\mathcal{L}$ whose members belong to $W$. Thus \eqref{eq10.30} is proved.

As in the proof of (1), we obtain the first inclusion \eqref{eq10.31} from \eqref{eq10.30}. Similarly, the second inclusion in \eqref{eq10.31} can be deduced from \eqref{eq10.30} and \eqref{eqP}. Then \eqref{eq10.32} follows from \eqref{eq10.30} and \eqref{eq10.31} by using the time derivative of \eqref{1.2'}.

Now by the Sobolev embedding, $\partial_x \phi\in L^\infty_t(H^1(\Omega))\subset L^\infty_t(L^6(\Omega))$. By differentiating \eqref{1.2'} and \eqref{eqP} with respect to $x$ and applying Lemma \ref{lem3.1}, we have $u_x\in L^\infty_t(L^6(\Omega))$. By Lemma \ref{remark-imbedding}, $\nabla \phi\in L^\infty_t(L^6(\Omega))$, and thus $\partial_x u\cdot \nabla \phi\in L^\infty_t(L^3(\Omega))$. Now we differentiate \eqref{1.1'}-\eqref{1.3'} in $x$ to deduce 
\begin{equation}
    \label{eq11.34}
\frac{\partial}{\partial t}\frac{\partial \phi}{\partial x}
-\operatorname{div}( D \nabla \frac{\partial\phi}{\partial x})=-\frac{\partial\mathbf{u}}{\partial x} \cdot \nabla \phi-\mathbf{u} \cdot \nabla \frac{\partial\phi}{\partial x}, \quad \frac{\partial \phi}{\partial x}\big|_{z=0, -H}=0, \quad \frac{\partial\phi}{\partial x}\big|_{t=0}=\frac{\partial \phi_0}{\partial x}\in H^1.
\end{equation} 
Since $\mathbf{u}\in L^\infty_{t,x}$, $\mathbf{u} \cdot \nabla \frac{\partial\phi}{\partial x}\in L^2_{t,x}$, $\partial_x \mathbf{u}\cdot \nabla \phi\in L^\infty_t(L^3(\Omega))$, $\partial_x \phi_0\in H^1(\Omega)$, we can rewrite the equation \eqref{eq11.34} into a non-divergence form equation with the right-hand side in $L^2_{t,x}$ and initial data in $H^1$. Then by the $W^{1,2}_2$ estimate for parabolic equations with leading coefficients depending only on $z$, we get $\partial_x \phi\in L^2_t(W)$, which implies that $\partial_x \nabla \phi\Big|_{\Omega_j}\in L_t^2(H^1(\Omega_j))$ for each $j$. This together with $\nabla \phi\Big|_{\Omega_j}\in L_t^2(H^1(\Omega_j))$ gives $\nabla \phi\in L_t^2(L^\infty(\Omega))$ in view of Lemma \ref{lem6.1}. 
This ends the proof of (2).

(3) 
Using \eqref{eq11.34}, we repeat the same argument as those for the proof of (1) and (2). The equations \eqref{eq11.34} are similar to the original ones modulo an additional term $\partial_x\mathbf{u} \cdot \nabla \phi$, which is of lower order and in $L^\infty(0,T;L^3(\Omega))$ because $\frac{\partial\mathbf{u}}{\partial x},\nabla \phi\in L^\infty_t(L^6(\Omega))$ according to (2). Using the fact that $\partial_t u\partial_x \phi$ and $\partial_x u\partial_t \phi$ are in $L^2(0,T;L^3(\Omega))$, we deduce as before that
$$\frac{\partial^2 \phi}{\partial t\partial x} \in L^\infty_t(L^2(\Omega))\cap L^2_t(V),\quad \frac{\partial\phi}{\partial x}\in L^\infty(0,T; W).$$
 Hence,
$$\nabla\frac{\partial\phi}{\partial x}\big|_{\Omega_j}  \in L^\infty(0,T; H^1(\Omega_j)),\quad \forall j.$$
By Lemma \ref{lem6.1}, we conclude
$$\|\nabla\phi\|_{L^\infty(0,T;L^\infty(\Omega_j))} 
\le C, \quad \forall j.$$
Therefore, \eqref{eq10.33} is proved.
As in the proofs of (1) and (2), \eqref{eq10.34} and \eqref{eq10.35} follows from \eqref{eq10.33} together with \eqref{eqP} and \eqref{1.2'}.
Next we differentiate \eqref{eq11.34} again with respect to $x$ to get
$$
\frac{\partial}{\partial t}\partial_x^2\phi 
-\operatorname{div}( D \nabla \partial_x^2\phi)=-\partial_x^2\mathbf{u} \cdot \nabla \phi-\mathbf{u} \cdot \nabla \partial_x^2\phi-2\partial_x\mathbf{u} \cdot \nabla \partial_x\phi
$$
with the corresponding boundary condition and the initial condition $\partial_x^2 \phi_0\in H^1(\Omega)$. By the estimates above, it is easily seen that the right-hand side above is in $L^2_{t,x}$. Therefore, arguing as before, we get $\partial_x^2\phi \in L_t^2(W)$, which implies that $\partial_x^2 \nabla u\Big|_{\Omega_j}\in L_t^2(H^1(\Omega_j))$. This together with $\partial_x\nabla \phi\Big|_{\Omega_j}\in L_t^2(H^1(\Omega_j))$ gives $\partial_x\nabla \phi\in L_t^2(L^\infty(\Omega))$.
This completes the proof of (3). 

(4) We first consider the case when $k=2$. We differentiate the equations \eqref{2.1'}-\eqref{2.3'} in $x$ twice to deduce
$$ 
\frac{\partial}{\partial t}\frac{\partial^2 \phi}{\partial x^2}+\mathbf{u} \cdot \nabla \frac{\partial^2 \phi}{\partial x^2}
+\frac{\partial^2 \mathbf{u}}{\partial x^2} \cdot \nabla \phi
+2\frac{\partial \mathbf{u}}{\partial x} \cdot \nabla \frac{\partial\phi}{\partial x}
-\operatorname{div}( D \nabla \frac{\partial^2\phi}{\partial x^2})=0
$$
with the boundary and initial conditions
$$
\frac{\partial^2 \phi}{\partial x^2}\big|_{z=0, -H}=0, \quad \frac{\partial^2 \phi}{\partial x^2}\big|_{t=0}=\frac{\partial^2 \phi_0}{\partial x^2}\in W.
$$
Now we estimate the lower-order terms $\frac{\partial^2 \mathbf{u}}{\partial x^2} \cdot \nabla \phi$ and $\frac{\partial \mathbf{u}}{\partial x} \cdot \nabla \frac{\partial\phi}{\partial x}$, which belong to $L^\infty_t(L^6(\Omega))$ according to (3). Then we differentiate the above equation in $t$ and observe that the nonlinear terms
$$
\partial_t \mathbf{u}\frac{\partial^2 \phi}{\partial x^2},\quad
\partial_t \frac{\partial^2 \mathbf{u}}{\partial x^2}\phi,\quad
\frac{\partial^2 \mathbf{u}}{\partial x^2} \partial_t\phi,\quad
\partial_t\frac{\partial \mathbf{u}}{\partial x} \frac{\partial\phi}{\partial x},\quad
\frac{\partial \mathbf{u}}{\partial x} \partial_t\frac{\partial\phi}{\partial x}
$$
all belong to $L^2_t(L^2(\Omega))$ according to (2) and (3). By testing the equation of $\partial_t\partial_x^2 \phi$ by itself and integrating by parts as before, we deduce that 
$$
\partial_t\partial_x^2 \phi\in L^\infty_t(L^2(\Omega))\cap L^2_t(V).
$$
Now by using the elliptic regularity, we further obtain
$$
\partial_x^2 \phi\in L^\infty(0,T; W).
$$
 Hence,
$$\nabla\partial_x^2\phi\big|_{\Omega_j}  \in L^\infty(0,T; H^1(\Omega_j)),\quad \forall j.$$
By Lemma \ref{lem6.1}, we conclude
$$\|\nabla\partial_x \phi\|_{L^\infty(0,T;L^\infty(\Omega_j))} 
\le C, \quad \forall j.$$ 
Thus, \eqref{eq11.45} is proved.
As before, \eqref{eq11.46} and \eqref{eq11.47} follows from \eqref{eq11.45} together with \eqref{eqP} and \eqref{1.2'}.
Finally, we deduce \eqref{eq11.48} by using the same argument like before.
For general $k$, we use an induction argument. This completes the proof of (4).
\end{proof}

\subsection{Estimates on the diffuse interface model}
\begin{lemma}
 \label{lem-diffuse}
The solution to the diffuse interface model \eqref{2.1'}-\eqref{2.3'} enjoys the following properties:
\begin{itemize}
  \item For $\phi_0\in W$, the solution is uniformly bounded in space and time for bounded initial data, i.e., there exists a constant $C$, s.t.
  \begin{equation}
            \label{eq11.53}
	\|\phi^\veps\|_{L^\infty(0,T; L^\infty(\Omega))} \le C.
    \end{equation}
  \item For $\phi_0, \frac{\partial\phi_0}{\partial x}\in W$,  there exists $C>0$, independent of $\veps$, s.t.
  \be \label{eq11.54}
  \|\partial_x \phi^\veps\|_{L^\infty(0,T; L^\infty(\Omega))}
  \le C.\ee
  \item For $\phi_0, \frac{\partial\phi_0}{\partial x}, \frac{\partial^2\phi_0}{\partial x^2}\in W$,  there exists $C>0$, independent of $\veps$, s.t.
  \be \label{eq11.55}
  \|\frac{\partial^2\phi^\veps}{\partial x^2}\|_{L^\infty(0,T; L^\infty(\Omega))} \le C. \ee
\ignore{
  \item {\color{red} This point is not needed at this time} In the case when the sharp interface is located at $z=0$, and $\phi_0\in ???$, the pressure is uniformly bounded in $\veps$ with weight, i.e., $\exists C>0$, s.t.
  \begin{equation}
    \|z^k p^\veps\|_{L^\infty(0,T; W^{k+1,q}(\Omega))}\le C.
  \end{equation}
  }
\end{itemize}
\end{lemma}
\begin{proof}
    First, \eqref{eq11.53} follows from the maximum principle for parabolic equations. To prove \eqref{eq11.54}, we apply the Moser iteration to 
$$ 
\frac{\partial}{\partial t}\frac{\partial \phi^\veps}{\partial x}+\mathbf{u}^\veps \cdot \nabla \frac{\partial\phi^\veps}{\partial x}
+\frac{\partial\mathbf{u}^\veps}{\partial x} \cdot \nabla \phi^\veps
-\operatorname{div}( D^\veps \nabla \frac{\partial\phi^\veps}{\partial x})=0, \quad \frac{\partial \phi^\veps}{\partial x}\big|_{z=0, -H}=0, \quad \frac{\partial\phi^\veps}{\partial x}\big|_{t=0}=\frac{\partial \phi_0}{\partial x}\in W.$$
We test the equation above with $|\partial_x\phi^\veps|^{p-2}\partial_x\phi^\veps$, where $p\ge 2$. For the third term on the left-hand side, we use the boundedness of $\phi^\veps$ to get
\begin{eqnarray*}
&\Big|\int_\Omega \partial_x \mathbf{u}^\veps\cdot\nabla \phi^\veps |\partial_x\phi^\veps|^{p-2}\partial_x\phi^\veps\,dx\Big|
=(p-1)\Big|\int_\Omega \partial_x \mathbf{u}^\veps \phi^\veps |\partial_x\phi^\veps|^{p-2}\cdot\nabla\partial_x\phi^\veps\,dx\Big|    \\
&\le C(p-1)/p\Big|\int_\Omega |\partial_x \mathbf{u}^\veps||\nabla|\partial_x\phi^\veps|^{p/2}||\partial_x\phi^\veps|^{p/2-1}\,dx\Big|\\
&\le \frac{\min D_j}{2p}\int_\Omega |\nabla |\partial_x \phi|^{p/2}|^2+Cp\|\partial_x \mathbf{u}^\veps\|_{L^p}^2\|\partial_x \phi^\veps\|_{L^p}^{p-2}\\
&\le \frac{\min D_j}{2p}\int_\Omega |\nabla |\partial_x \phi|^{p/2}|^2+Cp^\beta\|\partial_x \phi^\veps\|_{L^p}^{p}
\end{eqnarray*}
for some constant $\beta\ge 1$, where the last inequality is due to Lemma \ref{lem3.1} and \eqref{2.2'}. With this and the divergence free condition on $\mathbf{u}$, we obtain
$$
\partial_t \frac 1 p\int_{\Omega}|\partial_x \phi^\varepsilon|^p\,dx
+\frac {\min D_j} p\int_{\Omega}|\nabla |\partial_x \phi^\varepsilon|^{p/2}|^2\,dx
\le Cp^\beta\|\partial_x \phi^\veps\|_{L^p_x}^{p}
$$
and thus
\begin{eqnarray*}
&\sup_{t\in [0,T]} \int_{\Omega}|\partial_x \phi^\varepsilon(t,x)|^p\,dx
+\min D_j\int_0^T\int_{\Omega}|\nabla |\partial_x \phi^\varepsilon|^{p/2}|^2\,dxdt\\
&\le Cp^{\beta+1}\|\partial_x \phi^\veps\|_{L^p(0,T;L^p(\Omega)}^{p}
+\|\partial_x \phi_0\|_{L^p(\Omega)}^p.    
\end{eqnarray*}
It then follows from the parabolic embedding 
$$
L^\infty(0,T;L^2)\cap L^2(0,T;H^1)\subset L^{2(d+2)/d}_{t,x},\quad
\|f\|_{L^{2(d+2)/d}_{t,x}}\le C\|f\|_{L^\infty(0,T;L^2)}+C\|f\|_{L^2(0,T;H^1)}
$$
that
$$
\|\partial_x \phi^\veps\|_{L^{p\gamma}(0,T;L^{p\gamma}(\Omega)}
\le Cp^{(\beta+1)/p}\|\partial_x \phi^\veps\|_{L^p(0,T;L^p(\Omega)}
+C|\Omega|^{1/p}\|\partial_x \phi_0\|_{L^\infty(\Omega)},
$$
where $\gamma=(d+2)/d>1$.
Then by the standard Moser iteration and using the boundedness of $\partial_x \phi_0$, we derive the boundedness of $\partial_x \phi^\veps$.

Similarly, to prove \eqref{eq11.55}, we apply the Moser iteration to
$$ 
\frac{\partial}{\partial t}\frac{\partial^2 \phi^\veps}{\partial x^2}+\mathbf{u}^\veps \cdot \nabla \frac{\partial^2\phi^\veps}{\partial x^2}
+\frac{\partial^2\mathbf{u}^\veps}{\partial x^2} \cdot \nabla \phi^\veps
+2\frac{\partial\mathbf{u}^\veps}{\partial x} \cdot \nabla \frac{\partial\phi^\veps}{\partial x}
-\operatorname{div}( D \nabla \frac{\partial^2\phi^\veps}{\partial x^2})=0
$$
with the corresponding boundary and initial conditions.
We test the equation above with $|\partial^2_x\phi^\veps|^{p-2}\partial^2_x\phi^\veps$. The third term can be bounded as before. For the fourth term on the left-hand side, we again integrate by parts and use the boundedness of $\partial_x \phi^\veps$ and the $L^p$ bound of $\partial_x \mathbf{u}^\veps$ as in the previous step.
\end{proof}

\begin{remark}
    Notice that the proof above covers the case of less regular initial data with slightly weaker results. In particular, if  $\frac{\partial\phi_0}{\partial x}\in V$, then $\forall k\in \mathbb{N}^+$,  there exists $C_k>0$, independent of $\veps$, s.t.
$
  \|\partial_x \phi^\veps\|_{L^\infty(0,T; L^{2k}(\Omega))}  \le C_k$.
 If $\frac{\partial\phi_0}{\partial x}, \frac{\partial^2\phi_0}{\partial x^2}\in V$, then $\forall  k\in \mathbb{N}^+$, there exists $C_k>0$, independent of $\veps$, s.t.
$ 
  \|\frac{\partial^2\phi^\veps}{\partial x^2}\|_{L^\infty(0,T; L^{2k}(\Omega))} \le C_k. 
  $
  \end{remark}
\ignore{
For the proof of \eqref{eq11.54'}, we test the $x$-derivative of equation \eqref{2.3'} by $(\partial_x\phi^\veps)^{2k-1}$. We have two nonlinear terms to deal with,  $\mathbf{u}^\veps \cdot \nabla \partial_x\phi^\veps
$, and $\partial_x\mathbf{u}^\veps\cdot \nabla \phi^\veps$. The first term yields zero after an integration by parts. For the second term, we have the following.
\begin{align*}
&\Big|\int_\Omega \partial_x \mathbf{u}^\veps\cdot\nabla \phi^\veps (\partial_x\phi^\veps)^{2k-1}\,dx\Big|
=(2k-1)\Big|\int_\Omega \partial_x \mathbf{u}^\veps \phi^\veps (\partial_x\phi^\veps)^{2k-2}\cdot\nabla\partial_x\phi^\veps\,dx\Big|    
\\
&\le C (2k-1)\|\partial_x \mathbf{u}^\veps\|_{L^{2k}}\|\partial_x\phi^\veps\|^{k-1}_{L^{2k}} \|(\partial_x\phi^\veps)^{k-1}\nabla\partial_x \phi^\veps\|_{L^2}
\\
&\le \frac{2k-1}{2}\|\sqrt{D^\veps}(\partial_x\phi^\veps)^{k-1}\nabla\partial_x \phi^\veps\|^2_{L^2}+Ck\|\partial_x \mathbf{u}^\veps\|_{L^{2k}}^2\|\partial_x \phi^\veps\|_{L^{2k}}^{2k-2}
\\
&\le \frac{2k-1}{2}\|\sqrt{D^\veps}(\partial_x\phi^\veps)^{k-1}\nabla\partial_x \phi^\veps\|^2_{L^2}+Ck C(k)\|\partial_x \phi^\veps\|_{L^{2k}}^{2k},
\end{align*}
where $C(k)$ is the constant for the $W^{1,2k}$ elliptic estimates associated with \eqref{2.2'}.
Therefore,
$$\frac{1}{k}\frac{d}{dt}\|\partial_x\phi^\veps\|_{L^{2k}}^{2k} +(2k-1)\|\sqrt{D^\veps}(\partial_x\phi^\veps)^{k-1}\nabla\partial_x \phi^\veps\|^2_{L^2}\le Ck C(k)\|\partial_x \phi^\veps\|_{L^{2k}}^{2k}. $$
The desired result follows from Gronwall.

Similarly, to prove \eqref{eq11.55'}, we perform energy estimates to
$$ 
\frac{\partial}{\partial t}\frac{\partial^2 \phi^\veps}{\partial x^2}+\mathbf{u}^\veps \cdot \nabla \frac{\partial^2\phi^\veps}{\partial x^2}
+\frac{\partial^2\mathbf{u}^\veps}{\partial x^2} \cdot \nabla \phi^\veps
+2\frac{\partial\mathbf{u}^\veps}{\partial x} \cdot \nabla \frac{\partial\phi^\veps}{\partial x}
-\operatorname{div}( D^\veps \nabla \frac{\partial^2\phi^\veps}{\partial x^2})=0
$$
with the corresponding boundary and initial conditions.
We test the equation above with $(\partial^2_x\phi^\veps)^{2k-1}$. The third term can be bounded as before. For the fourth term on the left-hand side, we again integrate by parts and use the $L^{4k}$ bound on $\partial_x \phi^\veps$ and  $\partial_x \mathbf{u}^\veps$ as in the previous step. More specifically, we have,
\begin{align*}
&\Big|\int_{\Omega} \partial_x\mathbf{u}^\veps\cdot\nabla\partial_x\phi^\veps (\partial^2_x\phi^\veps)^{2k-1}\,d\mathbf{x}\Big|
= (2k-1) \Big|\int_{\Omega} \partial_x\mathbf{u}^\veps\partial_x\phi^\veps (\partial^2_x\phi^\veps)^{2k-2}\nabla\partial^2_x\phi^\veps\,d\mathbf{x}\Big|
\\
& \le (2k-1) \|\partial_x\mathbf{u}^\veps\|_{L^{4k}}\|\partial_x\phi^\veps\|_{L^{4k}}
  \|\partial^2_x\phi^\veps\|_{L^{2k}}^{k-1}\|(\partial^2_x\phi^\veps)^{k-1}\nabla\partial^2_x\phi^\veps\|_{L^2}
  \\
  &\le \frac{2k-1}{4}\|\sqrt{D^\veps}(\partial^2_x\phi^\veps)^{k-1}\nabla\partial^2_x\phi^\veps\|^2_{L^2}
  + Ck  \|\partial_x\mathbf{u}^\veps\|^2_{L^{4k}}\|\partial_x\phi^\veps\|^2_{L^{4k}}
  \|\partial^2_x\phi^\veps\|_{L^{2k}}^{2k-2}.
\end{align*}
The desired \eqref{eq11.55'} then follows from \eqref{eq11.54'} and straightforward Gronwall.
}

There is a weak formulation of the diffuse material interface model similar to the weak formulation of the sharp material interface model \eqref{def-eq} as follows:
\begin{equation}\label{def-eq-d}
\left(\phi^\veps(t_2), \psi\right)-\left(\phi^\veps(t_1), \psi\right)+\int_{t_1}^{t_2}A^\veps( \phi^\veps, \psi) dt+\int_{t_1}^{t_2}B(\mathbf{u}^\veps,  \phi^\veps,\psi) d t=0,
\end{equation}
where
\begin{equation}
A^\veps(\phi,\psi )=  \left( b D^\veps\nabla \phi, \nabla \psi\right), \forall \phi, \psi \in V , \quad
B(\mathbf{u}, \phi,\psi )= \int_{\Omega}\left(\mathbf{u}\cdot \nabla \phi \right)\psi d\bx, \forall \phi, \psi \in V. 
\end{equation}
For regular enough initial data, weak solution of the diffuse material interface model coincides with strong solution. 
We will utilize the weak formulation to investigate the convergence. 

 \section{Convergence}

 We now focus on the convergence of the solutions of the diffuse model to that of the sharp interface model as the width of the transition layer for the material properties approach zero.
 We first show the uniform convergence of the concentration and the pressure, as well as the $L^2$ convergence of the velocity. 
 We then demonstrate that the tangential velocity is discontinuous across the sharp interfaces in the generic case. We construct an approximate velocity that is the sharp interface solution plus a boundary layer type near singular term. We rely on a novel anisotropic Sobolev imbedding, proved in the appendix,  to show that the approximate velocity is close to the velocity field of the diffuse model uniformly.
 \subsection{Convergence of the concentration \texorpdfstring{$\phi$}{} and the pressure \texorpdfstring{$P$}{}}
 Since no singularity is expected in the concentration, we can look at the difference of the concentrations directly.
 Let
 $$\delta\phi^\veps = \phi^\veps -\phi,\quad \delta p^\veps = p^\veps - p, \quad\delta\mathbf{u}^\veps =\mathbf{u}^\eps - \mathbf{u},\quad \delta K^\veps=K^\veps-K, \quad\delta D^\veps=D^\veps-D.$$
 Taking the difference between \eqref{1.1'}-\eqref{1.3'} and \eqref{2.1'}-\eqref{2.3'}, we deduce the following equation for the differences
 \begin{eqnarray}
\operatorname{div}(\delta\mathbf{u}^\veps)=0,  \label{d.1}\\
\delta\mathbf{u}^\veps=-K^\veps\left(\nabla \delta P^\veps+ \delta\phi^\veps \mathbf{e}_z \right)-\delta K^\veps\left(\nabla P+ \phi \mathbf{e}_z \right),\label{d.2} \\
 \frac{\partial \delta\phi^\veps}{\partial t}+\mathbf{u}^\veps \cdot \nabla \delta\phi^\veps +\delta\mathbf{u}^\veps \cdot \nabla \phi
 		-\operatorname{div}( D^\veps \nabla \delta\phi^\veps) -\operatorname{div}( \delta D^\veps \nabla \phi)=0.\label{d.3}
\end{eqnarray}
These equations should be understood in the weak sense utilizing the weak formulations \eqref{def-eq} and \eqref{def-eq-d}.
The pressure and velocity differences may be estimated via the concentration difference and permeability difference using \eqref{d.1}-\eqref{d.2} together with standard elliptic regularity:
\begin{eqnarray*}
\| \sqrt{K^\veps} \nabla\delta P^\veps\|_{L^2}
&\le& C(\|\sqrt{K^\veps} \delta\phi^\veps \|_{L^2} + \|\delta K^\veps\left(\nabla P+ \phi \mathbf{e}_z \right)\|_{L^2})
\\
&\le& C \|\delta\phi^\veps \|_{L^2} + C\|\delta K^\veps\|_{L^2}\left(\|\nabla P\|_{L^\infty}+ \|\phi \|_{L^\infty}\right)\\
&\le& C \|\delta\phi^\veps \|_{L^2} + C\veps^{1/2}
\end{eqnarray*}
and
\begin{eqnarray*}
\| \delta\mathbf{u}^\veps\|_{L^2}
&\le& \|K^\veps\|_{L^\infty} (\|\nabla\delta P^\veps\|_{L^2}+ \|\delta\phi^\veps\|_{L^2})
		+ \|\delta K^\veps\|_{L^2}\left(\|\nabla P\|_{L^\infty}+ \|\phi \|_{L^\infty}\right)
		\\
&\le & C \|\delta\phi^\veps \|_{L^2} + C\|\delta K^\veps\|_{L^2}\left(\|\nabla P\|_{L^\infty}+ \|\phi \|_{L^\infty}\right)
\\
&\le & C \|\delta\phi^\veps \|_{L^2} + C\veps^\frac12,
\end{eqnarray*}
where we have used  the uniform bound on $\phi$ and $\nabla P$ and the assumption that the permeability difference is supported in a narrow band of width $\veps$.

For the estimate on the difference of concentrations, we multiply \eqref{d.3} by $\delta\phi^\veps$ and integrate over the domain, and we deduce
\begin{eqnarray*}
 \frac12\frac{d}{d t}\| \delta\phi^\veps\|^2_{L^2} + \|\sqrt{D^\veps} \nabla \delta\phi^\veps\|^2_{L^2}
 &\le&
 	\|\delta\mathbf{u}^\veps\|_{L^2} \| \nabla \delta\phi^\veps\|_{L^2}\|\phi\|_{L^\infty}
 		+\| \delta D^\veps\|_{L^2}\| \nabla \phi\|_{L^\infty} \| \nabla \delta\phi^\veps\|_{L^2}
\\
&\le&
 	 \frac12\|\sqrt{D^\veps} \nabla \delta\phi^\veps\|^2_{L^2}
 	 	+ C \|\delta\mathbf{u}^\veps\|^2_{L^2}\|\phi\|^2_{L^\infty}
 	 	+ C \| \delta D^\veps\|^2_{L^2}\| \nabla \phi\|^2_{L^\infty}.
\end{eqnarray*}
This implies, with the uniform estimate on $\nabla\phi$ and the assumption that the diffusivity difference is supported in a narrow band of width $\veps$, the first part of the following convergence result.
\begin{theorem}
\label{thm-conv-phi}
 For regular enough initial data, the solution to the diffuse interface model converges to the corresponding solution of the sharp interface model at vanishing material interface width. More specifically, we have, 
\begin{enumerate}
 \item If $\phi_0, \partial_x\phi_0\in W$, 
$\exists C>0$, independent of $\veps$, such that 
\begin{eqnarray}
                \label{eq4.4}
\|\delta\phi^\veps\|_{L^\infty(0, T; L^2)} + \|\frac{\partial}{\partial x}\delta\phi^\veps\|_{L^\infty(0, T; L^2)}&\le & C\veps^\frac12,
\\
                \label{eq4.5}
\|\nabla\delta\phi^\veps\|_{L^2(0, T; L^2)} + \|\nabla\frac{\partial}{\partial x}\delta\phi^\veps\|_{L^2(0, T; L^2)} &\le & C\veps^\frac12,
\\
                \label{eq4.6}
\|\nabla\delta P^\veps\|_{L^\infty(0, T; L^2)} + \|\nabla\frac{\partial}{\partial x}\delta P^\veps\|_{L^\infty(0, T; L^2)} &\le & C\veps^\frac12,
\\
                \label{eq4.7}
\|\delta\mathbf{u}^\veps\|_{L^\infty(0, T; L^2)} + \|\frac{\partial}{\partial x}\delta\mathbf{u}^\veps\|_{L^\infty(0, T; L^2)}&\le & C\veps^\frac12,
\\
    \label{u-conv}
\|\delta\phi^\veps\|_{L^2_t(L^\infty)} + \|\delta P^\veps\|_{L^\infty_t(L^\infty)} &\le& C\veps^\frac12. 
\end{eqnarray}
Moreover, we can take the vanishing material interface width limit in \eqref{def-eq-d} to recover \eqref{def-eq} together with the interfacial boundary conditions \eqref{IBC}.

  \item If $\phi_0, , \partial_x\phi_0 , \partial^2_x\phi_0 \in W$, then 
$\exists C>0$, independent of $\veps$, such that
\begin{eqnarray}
 \|\frac{\partial^2}{\partial x^2}\delta\phi^\veps\|_{L^\infty(0, T; L^2)} &\le & C\veps^\frac12,
\\
            \label{eq4.10}
\|\nabla\frac{\partial^2}{\partial x^2}\delta\phi^\veps\|_{L^2(0, T; L^2)} &\le & C\veps^\frac12,
\\
\|\nabla\frac{\partial^2}{\partial x^2}\delta P^\veps\|_{L^\infty(0, T; L^2)} &\le & C\veps^\frac12, \label{P-dxx}
\\
\|\frac{\partial^2}{\partial x^2}\delta\mathbf{u}^\veps\|_{L^\infty(0, T; L^2)} &\le & C\veps^\frac12,
\\
\label{u-conv-x}
\|\partial_x\delta\phi^\veps\|_{L^2_t(L^\infty)} + \|\partial_x\delta P^\veps\|_{L^\infty_t(L^\infty)} &\le& C\veps^\frac12. 
\end{eqnarray}
\end{enumerate}
\end{theorem}

\begin{proof}
(1) We have already bounded the first terms in  \eqref{eq4.4} - \eqref{eq4.7}. The interfacial boundary conditions \eqref{IBC} are embedded in the solution space and weak formulation of the models.
To bound the second terms, we differentiate \eqref{d.1} - \eqref{d.3} in $x$ to get
 \begin{eqnarray}
&\operatorname{div}(\partial_x \delta\mathbf{u}^\veps)=0,  \label{d.1b}\\
&\partial_x\delta\mathbf{u}^\veps=-K^\veps\left(\nabla\delta\partial_x P^\veps+ \partial_x\delta\phi^\veps \mathbf{e}_z \right)-\delta K^\veps\left(\nabla \partial_x P+ \partial_x \phi \mathbf{e}_z \right),\label{d.2b} \\
 &\frac{\partial \partial_x\delta\phi^\veps}{\partial t}+\mathbf{u}^\veps \cdot \nabla \partial_x \delta\phi^\veps
 +\partial_x\mathbf{u}^\veps \cdot \nabla  \delta\phi^\veps 
 +\partial_x\delta\mathbf{u}^\veps \cdot \nabla \phi\notag\\
&\quad  +\delta\mathbf{u}^\veps \cdot \partial_x\nabla \phi
 		-\operatorname{div}( D^\veps \nabla \partial_x\delta\phi^\veps) -\operatorname{div}( \delta D^\veps \partial_x\nabla \phi)=0.\label{d.3b}
\end{eqnarray}
As before, we have
\begin{eqnarray*}
\| \sqrt{K^\veps} \nabla\partial_x \delta P^\veps\|_{L^2}
&\le& C(\|\sqrt{K^\veps} \partial_x \delta\phi^\veps \|_{L^2} + \|\delta K^\veps\left(\nabla\partial_x P+ \partial_x\phi \mathbf{e}_z \right)\|_{L^2})
\\
&\le& C \|\partial_x\delta\phi^\veps \|_{L^2} + C\|\delta K^\veps\|_{L^2}\left(\|\nabla \partial_x P\|_{L^\infty}+ \|\partial_x\phi \|_{L^\infty}\right)\\
&\le& C \|\partial_x\delta\phi^\veps \|_{L^2} + C\veps^{1/2}.
\end{eqnarray*}
\begin{eqnarray*}
\| \partial_x\delta\mathbf{u}^\veps\|_{L^2}
&\le& \|K^\veps\|_{L^\infty} (\|\nabla\delta \partial_x P^\veps\|_{L^2}+ \|\partial_x \delta\phi^\veps\|_{L^2})
		+ \|\delta K^\veps\|_{L^2}\left(\|\nabla\partial_x P\|_{L^\infty}+ \|\partial_x \phi \|_{L^\infty}\right)
		\\
&\le & C \|\partial_x \delta\phi^\veps \|_{L^2} + C\|\delta K^\veps\|_{L^2}\left(\|\nabla\partial_x P\|_{L^\infty}+ \|\partial_x\phi \|_{L^\infty}\right)
\\
&\le & C \|\partial_x\delta\phi^\veps \|_{L^2} + C\veps^\frac12,
\end{eqnarray*}
where we have used  the boundedness on $\partial_x\phi$ and $\nabla\partial_x P$  in Proposition \ref{prop3.3} (3) and the 
assumption that the permeability difference is supported in a narrow band of width $\veps$. 
Similarly, we deduce from the equation of $\partial_x\delta \phi^\veps$ that
\begin{eqnarray*}
 &\frac12\frac{d}{d t}\| \partial_x\delta\phi^\veps\|^2_{L^2} + \|\sqrt{D^\veps} \nabla \partial_x \delta\phi^\veps\|^2_{L^2}\\
 &\le \|\partial_x \mathbf{u}^\veps\|_{L^3} \|\delta \phi^\veps\|_{L^6}\|\nabla \partial_x \delta\phi^\veps\|_{L^2}+
 	\|\partial_x \delta\mathbf{u}^\veps\|_{L^2} \|\phi\|_{L^\infty} \| \nabla \partial_x \delta\phi^\veps\|_{L^2}
    +\|\delta\mathbf{u}^\veps\|_{L^2} \|\partial_x \phi\|_{L^\infty}\| \nabla \partial_x \delta\phi^\veps\|_{L^2}\\
 		&\quad +\| \delta D^\veps\|_{L^2}\|\partial_x \nabla \phi\|_{L^\infty} \| \nabla \partial_x \delta\phi^\veps\|_{L^2}
\\
&\le
 	 \frac12\|\sqrt{D^\veps} \nabla \partial_x\delta\phi^\veps\|^2_{L^2}
 	 +C\|\partial_x \mathbf{u}^\veps\|^2_{L^3} \|\delta \phi^\veps\|^2_{L^6}	+ C \|\partial_x\delta\mathbf{u}^\veps\|^2_{L^2}\|\phi\|^2_{L^\infty}+C\|\delta\mathbf{u}^\veps\|^2_{L^2} \|\partial_x \phi\|^2_{L^\infty}\\
 	 &\quad	+ C \| \delta D^\veps\|^2_{L^2}\|\partial_x \nabla \phi\|^2_{L^\infty}\\
&\le
 	 \frac12\|\sqrt{D^\veps} \nabla \partial_x\delta\phi^\veps\|^2_{L^2}
 	 +C\|\partial_x \phi^\veps\|^2_{L^3} \| \delta \phi^\veps\|^2_{H^1}	+ C (\|\partial_x\delta\phi^\veps \|_{L^2}^2 + \veps)
     +C\veps \|\partial_x \nabla \phi\|^2_{L^\infty},
\end{eqnarray*}
where we used the boundedness of $\|\phi\|_{L^\infty}$ and $\|\partial_x \phi\|_{L^\infty}$ guaranteed by Proposition \ref{prop3.3} (3). Recalling that $\partial_x \nabla \phi\in L^2_t(L^\infty(\Omega))$, 
$\|\partial_x \phi^\veps\|_{L^\infty_t(L^3)}\le C$ (cf. \eqref{eq11.54}), and $\| \delta\phi^\veps\|_{L^2_t(H^1)}\le C\veps^{1/2}$, we can conclude by applying the Gronwall inequality.

The desired uniform in space convergence of the concentration and pressure \eqref{u-conv} follows from \eqref{eq4.5}-\eqref{eq4.6} and lemma \ref{lem6.1}.

It is easy to see that we can take the vanishing $\veps$ limit in \eqref{def-eq-d} to recover the sharp interface model \eqref{def-eq} together with the interfacial boundary conditions.

(2) As in the previous step, using the boundedness of $\partial_x^2\phi$ and $\nabla \partial_x^2 P$ in Proposition \ref{prop3.3} (4) with $k=2$, we deduce that 
$$
\| \sqrt{K^\veps} \nabla\partial_x^2 \delta P^\veps\|_{L^2}
\le C \|\partial_x^2 \delta\phi^\veps \|_{L^2} + C\veps^{1/2},
$$
and
$$
\| \partial_x^2\delta\mathbf{u}^\veps\|_{L^2}
\le  C \|\partial_x^2\delta\phi^\veps \|_{L^2} + C\veps^\frac12.
$$
Note that $\partial_x^2\delta\phi^\veps$ satisfies
\begin{eqnarray*}
&\frac{\partial \partial_x^2\delta\phi^\veps}{\partial t}+\mathbf{u}^\veps \cdot \nabla \partial_x^2 \delta\phi^\veps
 +\partial_x^2\mathbf{u}^\veps \cdot \nabla  \delta\phi^\veps 
 +2\partial_x\mathbf{u}^\veps \cdot \nabla \partial_x \delta\phi^\veps
  +\partial_x^2\delta\mathbf{u}^\veps \cdot \nabla \phi\\
&\quad  +\delta\mathbf{u}^\veps \cdot \partial_x^2\nabla \phi
+2\partial_x \delta\mathbf{u}^\veps \cdot \partial_x\nabla \phi
 		-\operatorname{div}( D^\veps \nabla \partial_x^2\delta\phi^\veps) -\operatorname{div}( \delta D^\veps \partial_x^2\nabla \phi)=0.\label{d.3b}    
\end{eqnarray*}
By a similar reasoning, we get
\begin{eqnarray*}
 &\frac12\frac{d}{d t}\| \partial_x^2\delta\phi^\veps\|^2_{L^2} + \|\sqrt{D^\veps} \nabla \partial_x^2 \delta\phi^\veps\|^2_{L^2}\\
&\le
 	 \frac12\|\sqrt{D^\veps} \nabla \partial_x^2\delta\phi^\veps\|^2_{L^2}
 	 +C\|\partial_x^2 \phi^\veps\|^2_{L^3} \| \delta \phi^\veps\|^2_{H^1}\\
&\quad     +C\|\partial_x \phi^\veps\|^2_{L^3} \| \partial_x \delta \phi^\veps\|^2_{H^1}
     + C (\|\partial_x^2\delta\phi^\veps \|_{L^2}^2 + \veps)
     +C\veps\|\partial^2_x \nabla \phi\|^2_{L^\infty}.
\end{eqnarray*}
Recall that $\partial^2_x\nabla\phi\in L^2_t(L^\infty(\Omega))$ is guaranteed by Prop. 3.3 (4) with $k=2$. By the Gronwall inequality and the estimates from the previous steps, we derive the first four estimates. The desired uniform in space convergence of the$x$-derivative of the pressure \eqref{u-conv-x} follows from \eqref{eq4.5}, \eqref{eq4.6}, \eqref{eq4.10}, \eqref{P-dxx}, and lemma \ref{lem6.1}. This ends the proof.
\end{proof}

 \subsection{Boundary layer behavior of the velocity \texorpdfstring{$\mathbf{u}$}{}}
Although both $\phi$ and $P$ are continuous across the interfaces, the permeability $K$ is assumed to be discontinuous across the interfaces. Therefore, the tangential velocity, determined via \eqref{1.2'}, is generally discontinuous. In fact the discontinuity should satisfy the following jump condition
\begin{equation}
\label{v-jump}
 \llbracket u_x\rrbracket = -\frac{\partial P}{\partial x}\llbracket K\rrbracket  \mbox{ at}\ z=z_j, j=1,\cdots,\ell-1,
\end{equation}
although the vertical velocity must be continuous across the interfaces due to \eqref{1.2'}, i.e.,
\begin{equation}
\phi \llbracket K\rrbracket = -\llbracket K\frac{\partial P}{\partial z}\rrbracket  \mbox{ at}\ z=z_j, j=1,\cdots,\ell-1.
\end{equation}
Therefore, the velocity can not converge uniformly in space, or in any space with trace. This is why the sharp material interface limit is a singular limit from the mathematical point of view.
Moreover, the velocity converges uniformly away from the interface, as we shall illustrate below.
Hence, there must exist a transition layer in the velocity in the generic case. 

Our goal in this subsection is to construct an approximate velocty, $\tilde{\mathbf{u}}^\veps$, which is close to the diffusive solution $\mathbf{u}^\veps$ in $H^\alpha, \forall \alpha\in (\frac12, 1)$.  In addition, we demonstrate that this approximate velocity is the sum of the sharp interface velocity plus a boundary layer type transition layer with $H^\alpha$ norm. Notice that $H^\alpha, \alpha\in (\frac12, 1)$ is strong enough to have a trace on the interfaces. Thus, the $H^\alpha$ norm of the boundary layer type structure is not small at vanishing $\veps$. This is in contrast to the $L^2$ norm.

The vertical velocity is continuous for both the sharp and diffuse interface models. 
The pressure is also continuous for both models. However, the vertical derivative of the pressure in the sharp interface model is in general discontinuous \cite{mckibbin1980jfm, CNW2024}.  
Therefore, we proceed with the construction of a transition layer by an expansion in $P^\veps$ and $K^\veps$.
This approach is different from the classical boundary layer-type approach since we do not expect disparity in the leading order in $P$.

For simplicity, we consider a vertically translated non-dimensional domain, and assume that we have two layers only with the interface located at $z=0$ and the transition layer in the permeability and diffusivity of width $\veps$. The permeabilities in the top and bottom region are denoted $K_{\pm}$, while the diffusivities in the two regions are denoted $D_{\pm}$ with the plus sign denoting the upper region.

We first observe that the velocity converge uniformly in space away from the interface. Indeed, let
$\Omega_{+\veps}=\Omega\cap \{z\ge \veps\},  \Omega_{-\veps}=\Omega\cap \{z\le \veps\}$, we have
$$\Delta P^\veps = -\frac{\partial \phi^\veps}{\partial z} \in H^1(\Omega_{\pm\veps}).$$
Hence,
$$\|P^\veps\|_{H^3(\Omega_{\pm\veps})} \le C(\|\frac{\partial \phi^\veps}{\partial z}\|_{H^1}+\|P^\veps\|_{H^1})\le C.$$
Consequently,
$$\|\mathbf{u}^\veps \|_{H^2(\Omega_{\pm\veps})} \le C.$$
Similarly, we also have
$$\mathbf{u} \in H^2(\Omega_{\pm}).$$
Therefore, thanks to \eqref{eq4.7} and interpolation,
\begin{equation}
\label{v-uniform-away}
\|\delta\mathbf{u}^\veps\|_{L^\infty(\Omega_{\pm\veps} )}\le C \|\delta\mathbf{u}^\veps\|^\frac12_{L^2(\Omega_{\pm\veps} )} \|\delta\mathbf{u}^\veps\|^\frac12_{H^2(\Omega_{\pm\veps} )} \le C\veps^\frac14.
\end{equation}
This demonstrates the uniform convergence of the velocity away from the interface.
Likewise, the $H^2$ norm of the pressure converges uniformly away from the interface.

\ignore{
To deal with the transition within the interface, we introduce the stretched variable 
\begin{equation}
\tilde{z} = \frac{z}{\veps}
\end{equation}
We also denote
\be
 \bar{K}=\frac{K_+ + K_-}{2}, \delta K = \frac{K_+ - K_-}{2}
 \ee
}

To construct the approximate velocity, we introduce the stream functions $\psi^\veps, \psi$ for the diffuse and sharp material interface models, respectively. Recall that the velocity is determined by the streamfunction through $\mathbf{u}^\veps=\nabla^\perp\psi^\veps, \mathbf{u}=\nabla^\perp\psi$. We assume that $\psi^\veps\big|_{\text{top,bottom}}=0, \psi\big|_{\text{top,bottom}}=0$, consistent with the no-penetration velocity boundary condition \eqref{hBC} at the top and bottom.
Taking the curl of the velocity equations \eqref{1.2'} and \eqref{2.2'} we deduce
\begin{eqnarray}
 -\Delta\psi^\veps &=&\frac{\partial K^\veps}{\partial z}\frac{\partial P^\veps}{\partial x}-K^\veps\frac{\partial\phi^\veps}{\partial x},
\\
 -\Delta\psi &=&\frac{\partial K}{\partial z}\frac{\partial P}{\partial x}-K\frac{\partial\phi}{\partial x}.
\end{eqnarray}
We then introduce the following approximate streamfunction $\tilde{\psi}^\veps$ and velocity $\tilde{\mathbf{u}}^\veps$:
\be
\label{sf-approximate}
 -\Delta\tilde{\psi}^\veps =\frac{\partial K^\veps}{\partial z}\frac{\partial P }{\partial x}-K\frac{\partial\phi}{\partial x},
\quad \tilde{\psi}^\veps\big|_{\text{top,bottom}}=0, \quad \tilde{\mathbf{u}}^\veps=curl\ \tilde{\psi}^\veps.
\ee
$\tilde{\psi}^\veps$ can be solved with the help of an appropriate Green's function. 

{
The difference between $\tilde{\psi}^\veps$ and ${\psi}$ satisfies the following equation
\be
 -\Delta(\tilde{\psi}^\veps-\psi) =(\frac{\partial K^\veps}{\partial z}-\frac{\partial K}{\partial z})\frac{\partial P }{\partial x},
\quad \tilde{\psi}^\veps-\psi\big|_{\text{top,bottom}}=0.
\ee
The most singular part of the solution is associated with the right-hand side given by $(\frac{\partial K^\veps}{\partial z}-\frac{\partial K}{\partial z})(\frac{\partial  P}{\partial x}\big|_{z=0})$. 
This part of the solution is continuous. However, the associated velocity difference $\tilde{\mathbf{u}}^\veps-\mathbf{u}$ contains a jump discontinuity at the interface $z=0$, inherited from the jump discontinuity of the velocity field $\mathbf{u}$ of the sharp interface model.
In addition, $\tilde{\mathbf{u}}^\veps-\mathbf{u}$ is small away from the interface ($|z|\ge \veps$) as can be seen from the solution formula, or by triangle inequality and the smallness of ${\mathbf{u}}^\veps-\mathbf{u}$ and $\tilde{\mathbf{u}}^\veps-\mathbf{u}^\veps$ (to be proved below). Hence,
$\tilde{\mathbf{u}}^\veps-\mathbf{u}$ is of boundary layer type near $z=0$.
}

We now show that $\tilde{\mathbf{u}}^\veps$ is a good approximation of $\mathbf{u}^\eps$ uniformly in space.
Forming the difference of $\psi^\veps, \tilde{\psi}^\veps$, denoting the difference as $\delta\psi^\veps$, we deduce,
\be
\label{sf-diff}
-\Delta\delta\psi^\veps = \frac{\partial K^\veps}{\partial z}(\frac{\partial P^\veps}{\partial x}-\frac{\partial P }{\partial x})
	-(K^\veps-K)\frac{\partial\phi^\veps}{\partial x} -K \frac{\partial\delta\phi^\veps}{\partial x}.
\ee
Notice that the right-hand side is bounded uniformly in $\veps$ in the $L^2$ norm since
\begin{eqnarray*}
\|\frac{\partial K^\veps}{\partial z}(\frac{\partial P^\veps}{\partial x}-\frac{\partial P }{\partial x})\|_{L^\infty_t(L^2)}
&\le& \frac{C}{\veps}\veps^\frac12\|\frac{\partial P^\veps}{\partial x}-\frac{\partial P }{\partial x}\|_{L^\infty_t(L^\infty)} \le C,
\\
\|(K^\veps-K)\frac{\partial\phi^\veps}{\partial x} \|_{L^\infty_t(L^2)}
&\le& \| (K^\veps-K)\|_{L^2} \|\frac{\partial\phi^\veps}{\partial x}\|_{L^\infty_t(L^{\infty})} \le C_k\veps^\frac{1}{2}, 
\\
 \|K \frac{\partial\delta\phi^\veps}{\partial x}\|_{L^\infty_t(L^2)}
&\le& C\veps^\frac12.
\end{eqnarray*}
Here we have utilized the uniform convergence of the tangential derivative of the pressure, the uniform boundedness of the tangential derivative of the diffusive concentration, and the convergence of the tangential derivative of the concentration. Notice that the uniform convergence of the $x$ derivative of the pressure is a consequence of the anisotropic embedding Lemma \ref{lem6.1} and the convergence estimate \eqref{P-dxx}.
Consequently, we have, thanks to standard elliptic regularity,
\be
  \|\mathbf{u}^\veps-\tilde{\mathbf{u}}^\veps\|_{L^\infty_t(H^1(\Omega))}
\le \|\psi^\veps-\tilde{\psi}^\veps\|_{L^\infty_t(H^2(\Omega))} 
\le \|RHS\eqref{sf-diff}\|_{L^\infty_t(L^2)}
\le C.
\ee
Similarly, estimating the $H^{-1}$ norm of the right-hand side of \eqref{sf-diff}, we derive,
\be
\|\mathbf{u}^\veps-\tilde{\mathbf{u}}^\veps\|_{L^\infty_t(L^2(\Omega))}
\le \|\psi^\veps-\tilde{\psi}^\veps\|_{L^\infty_t(H^1(\Omega))} \le C\veps^\frac12.
\ee
Consequently,
\be
\|\mathbf{u}^\veps-\tilde{\mathbf{u}}^\veps\|_{L^\infty_t(H^\alpha(\Omega))} + \|\tilde{\mathbf{u}}^\veps-{\mathbf{u}}\|_{L^\infty_t(L^2(\Omega))} \le C\veps^{\frac12(1-\alpha)}, \quad \forall \alpha\in (\frac12, 1).
\ee

We now estimate the $x$ derivative of the velocity difference. For this purpose, we differentiate \eqref{sf-diff} in $x$, and we deduce
\be
\label{sf-diff-x}
-\Delta\frac{\partial}{\partial x}\delta\psi^\veps = \frac{\partial K^\veps}{\partial z}(\frac{\partial^2 P^\veps}{\partial x^2}-\frac{\partial^2 P }{\partial x^2})
	-(K^\veps-K)\frac{\partial^2\phi^\veps}{\partial x^2} -K \frac{\partial^2\delta\phi^\veps}{\partial x^2}.
\ee
Thanks to the convergence result on the concentration and pressure Theorem \ref{thm-conv-phi} and the a priori estimates on the diffuse model presented in lemma \ref{lem-diffuse}, we have
\begin{eqnarray*}
\|\frac{\partial^2 P^\veps}{\partial x^2}-\frac{\partial^2 P }{\partial x^2}\|_{L^\infty_t(L^q(\Omega))}&\le& C_q \veps^\frac12, \quad \forall q\in [1,\infty),
\\
\|\frac{\partial^2\phi^\veps}{\partial x^2}\|_{L^\infty_t(L^{\infty}(\Omega))} &\le& C,
\\
\|\frac{\partial^2\delta\phi^\veps}{\partial x^2}\|_{L^\infty_t(L^2(\Omega))}
&\le& C\veps^\frac12.
\end{eqnarray*}
This implies,
$$\| RHS \eqref{sf-diff-x}\|_{L^\infty_t(L^2)} \le C(q)\veps^{\frac 1 {q^*}-\frac 1 2}
+ C\veps^\frac12.
$$
Therefore,
\be
  \|\frac{\partial\mathbf{u}^\veps}{\partial x}-\frac{\partial\tilde{\mathbf{u}}^\veps}{\partial x}\|_{L^\infty_t(H^1(\Omega))} 
  = \|\frac{\partial}{\partial x}\delta\psi^\veps\|_{L^\infty_t(H^2(\Omega))}
  \le C\veps^{1/q^*-1/2},\quad \forall q^*\in (2,\infty),
\ee
where $q^*$ is the dual of $q$ determined via $\frac{1}{q}+\frac{1}{q^*}=\frac12$.\\
Similarly, we can easily deduce that the $H^{-1}$ norm of right-hand side of \eqref{sf-diff-x} satisfies
\be
  \|RHS\eqref{sf-diff-x}\|_{L^\infty_t(H^{-1}(\Omega))} \le C \veps^\frac{1}{2}.
\ee
This implies
\be
  \|\frac{\partial\mathbf{u}^\veps}{\partial x}-\frac{\partial\tilde{\mathbf{u}}^\veps}{\partial x}\|_{L^\infty_t(L^2(\Omega))} 
  = \|\frac{\partial}{\partial x}\delta\psi^\veps\|_{L^\infty_t(H^1(\Omega))}
  \le C\veps^{\frac{1}{2}}.
\ee
Interpolation leads to
\be
 \|\frac{\partial\mathbf{u}^\veps}{\partial x}-\frac{\partial\tilde{\mathbf{u}}^\veps}{\partial x}\|_{L^\infty_t(H^\alpha(\Omega))}  
 \le C \veps^{(1-\alpha)/2+\alpha(1/q^*-1/2)}, \quad \forall q^*\in (2,\infty).
\ee
Therefore, thanks to the anisotropic embedding \eqref{imbedding-anisotropic}, for any $q^*\in (2,\infty)$, $\alpha\in (\frac12, 1)$, and $\forall k\in\mathbb{N}^+$,
\begin{equation}
                \label{u-uniform}
\|\mathbf{u}^\veps - \tilde{\mathbf{u}}^\veps\|_{L^\infty_t(L^\infty(\Omega))} \le  C\veps^{(1-\alpha)/2+\alpha(1/q^*-1/2)} 
= C\veps^{\frac12 -\alpha + \frac{\alpha}{q^*}}. 
\end{equation}
Notice that for $\alpha\in (\frac12, 1)$, we can always choose $q^*>2$, s.t. $q^* < \frac{2\alpha}{2\alpha-1}$. This implies that
$$\frac12-\alpha + \frac{\alpha}{q^*} > 0,$$
This guarantees the uniform convergence in space and in time for the approximate velocity field.
By taking $\alpha\searrow 1/2$ and $q^* \searrow 2$, 
we get the convergence rate $\varepsilon^{\frac{1}{4}^-}$.

To summarise, we have shown 
\begin{theorem}
\label{them-conv-u}
 If $\phi_0, \frac{\partial}{\partial x}\phi_0, \frac{\partial^2}{\partial x^2}\phi_0\in W$, we have the following convergence results on the velocity field.
 \begin{enumerate}
  \item The velocity to the sharp interface model \eqref{1.1'}-\eqref{1.3'} is discontinuous in general. 
  The discontinuity must satisfy the jump condition \eqref{v-jump} across each interface.\\ 
    \item  The velocity field of the diffuse interface model \eqref{2.1'}-\eqref{2.3'}, which is continuous everywhere,  converges uniformly to the discontinuous velocity of the sharp interface model away from the interface in the sense that \eqref{v-uniform-away} holds. \\
    \item  The velocity of the diffuse interface model exhibits boundary layer type behavior near each interface. More specifically, the streamfunction of the diffuse interface model can be approximated  by the solution to \eqref{sf-approximate} under the $W^{1,\infty}$ norm in the sense of \eqref{u-uniform}. 
    Moreover, the velocity field of this approximate streamfunction exhibits explicit boundary layer type behavior near the interface $z=0$. 
  \end{enumerate}
\end{theorem}
\section{Conclusion and Remarks}

We have shown that the solutions to the diffuse material interface model for convection in layered porous media converge to those of the sharp material interface model, together with the appropriate interfacial boundary conditions. The limit is a singular one because of the existence of a narrow transition layer in the velocity field near each interface. Our result justifies the sharp material interface model together with the commonly adopted heuristic interfacial boundary conditions starting from the physically more realistic diffuse material interface model. 

The more general setting of smooth, curved interfaces can be treated analogously by invoking curvilinear coordinates, as long as the distance between different interfaces is bounded below by a positive constant. 

Our result can also be viewed as a structural stability result of the Darcy-Boussinesq system with respect to the diffusivity and permeability. Our result differs from the classical ones \cite{straughan2008} since we allow our physical parameters to be discontinuous.

While the result presented in this paper has provided an affirmative answer to the validity of the sharp material interface model on any finite time interval, the validity question of the model on infinite time interval is still open. This is related to important questions such as whether the long-time averaged Nusselt number \cite{doering1998jfm, wang2008pd} can be captured by the simplified sharp interface model. The mathematical investigation would involve the convergence of the global attractor and the invariant measures \cite{wang2008cpam}, among others. 
In addition, some applications involve one or more relatively thin layers with small permeability \cite{hewitt2022jfm, hewitt2020jfm}.
The simultaneous limit of vanishing thickness and the permeability of one or more layers is a challenging problem.

These and other physically important problems are the subject of subsequent work.


%
\section*{Appendix}
We present here a novel anisotropic embedding for functions with fractional regularity. This embedding result is crucial in establishing the explicit boundary layer type behavior of the velocity field.
\begin{lemma}
                    \label{lem6.1}
    Assume $u, \frac{\partial u}{\partial x}\in H^\alpha(\mathbb{R}^2),  \alpha>\frac12$. Then, $u\in L^\infty(\mathbb{R}^2)$. Moreover, $\exists C_\alpha>0$, independent of $u$, such that
    \be
    \label{imbedding-anisotropic}
    \|u\|_{L^\infty} \le C_\alpha (\|u\|_{H^\alpha} + \|\frac{\partial u}{\partial x}\|_{H^\alpha}).
    \ee
\end{lemma}
\begin{proof} 
    Let $\hat{u}(\xi)$ be the Fourier transform of $u$. Our assumption implies that
    $$\|(1+|\xi|^{2\alpha})(1+\xi_1^2)|\hat{u}|^2(\xi)\|_{L^1} \le C( \|u\|_{H^\alpha}^2 + \|\frac{\partial u}{\partial x}\|_{H^\alpha}^2)$$
    Recall that 
$$
|f(\bx)| \le \frac{1}{2\pi}(\int_{\mathbb{R}^2}|(1+|\xi|^{2\alpha})(1+\xi_1^2)|\hat{u}|^2(\xi)|d\xi)^\frac12
\times (\int_{\mathbb{R}^2}\frac{1}{(1+|\xi|^{2\alpha})(1+\xi_1^2)}d\xi)^\frac12.$$
By a change of variables, we have
\begin{eqnarray*}
   \int_{\mathbb{R}}\frac{1}{(1+|\xi|^{2\alpha})}d\xi_2
 \le \frac C {1+|\xi_1|^{2\alpha-1}}.
\end{eqnarray*}
Therefore,
$$
\int_{\mathbb{R}^2}\frac{1}{(1+|\xi|^{2\alpha})(1+\xi_1^2)}d\xi
\le C\int_{\mathbb{R}}\frac{1}{(1+\xi_1^2)^{\alpha+1/2}}d\xi_1<\infty.
$$


This completes the proof.
\end{proof}
\noindent{\bf Remark:} From the proof above, it is easily seen that the result still holds when $d=3$  provided that  we replace $\|\frac{\partial u}{\partial x}\|_{H^\alpha}$ with $ \|\frac{\partial u}{\partial x}\|_{H^\alpha} + \|\frac{\partial u}{\partial y}\|_{H^\alpha}$.

\noindent{\bf Remark:} For application in our channel geometry case, we could first extend in the vertical $z$ direction and then truncate in the horizontal $x$ direction to obtain a function on the whole plan with the $H^\alpha$ norms bounded by the original $H^\alpha$ norms plus the $L^2$ norm. Therefore, the embedding inequality \eqref{imbedding-anisotropic} remains valid.

\noindent{\bf Remark:} Anisotropic embeddings with integer $\alpha$ are classical. For instance, such an inequality with $\alpha=1$ was utilized in demonstrating the uniform in space validity of the Prandtl boundary layer theory for the Navier-Stokes equations with non-characteristic boundary conditions in \cite{TW2002}.
%

\bibliographystyle{amsalpha}
\bibliography{layer}

\end{document}